\documentclass[review]{elsarticle}

\usepackage{lineno,hyperref}
\usepackage{comment}
\usepackage{amsmath,amsxtra,amsthm,amssymb,latexsym, amscd}
\usepackage[mathscr]{eucal}
\usepackage{float}

\usepackage{lineno,hyperref}

\modulolinenumbers[5]

\numberwithin{equation}{section}
\newtheorem{theorem}{Theorem}[section]
\newtheorem{lemma}[theorem]{Lemma}

\newtheorem{corollary}[theorem]{Corollary}
\newtheorem{definition}[theorem]{Definition}
\newtheorem{remark}[theorem]{Remark}
\newtheorem{example}[theorem]{Example}

\allowdisplaybreaks

\begin{document}

\begin{frontmatter}

\title{Stability of a one-predator two-prey system governed by nonautonomous differential equations}



\author[Linh]{Linh Thi Hoai Nguyen}
\ead{nthlinh[at]ist.osaka-u.ac.jp}

\author[Quang]{Quang Hong Ta}
\ead{tahongquangosaka[at]gmail.com}

\author[Ton]{T\^{o}n Vi$\hat{\d{e}}$t T\d{a}\corref{mycorrespondingauthor}} 
\cortext[mycorrespondingauthor]{The last author was supported by JSPS KAKENHI Grant Number 20140047. }
\ead{tavietton[at]agr.kyushu-u.ac.jp}

\address[Linh]{Department of Information and Physical Sciences\\Graduate School of Information Science and Technology, Osaka University\\
1-5 Yamadaoka, Suita-shi, Osaka 565-0871, Japan
}

\address[Quang]{Fukuoka Japanese Language School\\
1-1-33 Hakataekihigashi, Hakata-ku,  Fukuoka 812-0013 Japan}

\address[Ton]{Center for Promotion of International Education and Research \\
Faculty of Agriculture, Kyushu University\\
6-10-1 Hakozaki, Higashi-ku, Fukuoka 812-8581, Japan}

\begin{abstract}
A non-periodic version of the one-predator two-prey system model presented in 
[L.T.H. Nguyen, Q.H. Ta, T.V. T\d{a}, Existence and stability of periodic solutions of a Lotka-Volterra system, SICE International Symposium on Control Systems,  
Tokyo, Japan, 712-4 (2015) 1-6] is considered. First, we prove existence of unique positive solutions to the model. Second, we show existence of an invariant set, which suggests the survival of all species in the system. 
On the other hand, we show that when the densities of two prey species are quite small, the predator falls into decay.
Third,   we explore global asymptotic stability of the system  by using the Lyapunov function method. Finally, some numerical examples are given to illustrate our results.
\end{abstract}


\end{frontmatter}

\section{Introduction}

Wildlife conservation is one of the most challenging problems in ecological and biological sciences. Knowing the mechanism of predator-prey  systems is one of crucial importance for maintaining the ecological integrity of  population and for preserving biodiversity.

The fundamental issue in these problems is to predict the variation of density of each species caused 
other species. 
Thus, we must first understand 
predator-prey interactions. 

In ecology, predator-prey interactions can be expressed by functional responses. 
Various functional responses have been proposed: 
\begin{itemize}
  \item Holling's type I functional response.  It is a linear function of the density of prey and is used  in  classical  Lotka-Volterra models.  It was the first kind of functional response described.
  \item Holling's type II functional response. It is a function of the density of prey, and is of  the form of a rectangular hyperbola (see \cite{Holling}).
  \item Ratio-dependent functional response. It is a function of the ratio of the density of predator to the density of prey  (see \cite{Arditi}). This function may be unbounded.
  \item Beddington-DeAngelis functional response. It is a bounded function of the ratio of the density of predator to the density of prey. This kind of functional response is  independently  presented by 
Beddington and DeAngelis et al. (see  \cite{Be,DGO}).
\end{itemize}

The latter type is of particular interest, as it takes into account both the densities of  predator and prey species. In fact, some ecologists showed that in many  situations, especially when predators have to search for food, functional responses should depend on  these densities (see \cite{APS,Do,JE,S_G} and references therein).  Furthermore, this type avoids some of  singular behaviors of ratio-dependent models at low prey density levels. In fact,   Beddington-DeAngelis functional responses have been used by many researchers (see  \cite{CC1,CC3,H1,TVT} and references therein).

We are interested in a one-predator two-prey system, where the two prey species are competitive with each other and are hunted by a predator species. Such a system may describe many biological  phenomena  in the real world. Let us give here two examples. (For more implication of three-species systems, we refer to 
\cite{Elettreby,Sharma,TVT2}  and references therein.)

\begin{example} 
Observe  a system of lion, buffalo, and  zebra species in the Serengeti ecosystem in Africa. On one hand, buffalo  and zebra species compete with each other for the food (plant). On the other hand, they are all attached by lion species. Is there a threshold for buffalo and zebra species for which if their density is smaller than the threshold, then the lion species heads towards extinction? In other words, can we find a quantitative relation
between the density of buffalo and zebra species and the dead rate of the lion species? The answer is affirmative (see Remark \ref{rm1}).
\end{example}

\begin{example} 
Consider mixed coniferous deciduous forests in the world (e.g., the mixed coniferous deciduous forest of white birch aspens fir and spruce trees in the Superior National Forest,  Minnesota, United States (\cite{Minnesota})). Coniferous and deciduous trees compete with each other for space. In fact, both of them  need space to get more light from the sun for their photosynthesis. Since we humans  exploit  the forests, humans can be assigned as a ``predator" species with respect to  both deciduous  and coniferous trees in the forest. This then makes a system of one predator and two prey species. 
\end{example}

In \cite{SICE}, we proposed a mathematical model for a one-predator two-prey population  governed by the  Beddington-DeAngelis functional response.
We solved the problem of periodicity in the system. However, the study only on periodicity seems to be incomplete. In fact, parameters in the model surely fluctuate due to the effect of internal and external factors. Some may be the competition between individuals of each species. Others may be the climate change.

In this paper, we consider a non-periodic version of the model presented in   \cite{SICE}. Our  model is performed by nonautonomous differential equations: 
\begin{equation}\label{E1}
\begin{cases}
\begin{aligned}
x_1'=&x_1\left[a_1(t)-b_{11}(t) x_1-b_{12}(t)x_2- \frac{\sigma_1(t)x_3}{\alpha(t)+\beta(t) x_1+\gamma(t) x_3}  \right ],\\
x_2'=&x_2\left[ a_2(t)-b_{21}(t) x_1-b_{22}(t)x_2- \frac{\sigma_2(t)x_2x_3}{\alpha(t)+\beta(t) x_2+\gamma(t) x_3}    \right ],\\
x_3'=&x_3\Big[-a_3(t)+ \frac{\rho_1(t)x_1}{\alpha(t)+\beta(t) x_1+\gamma(t) x_3}   \\
&\hspace{0.5cm} +\frac{\rho_2(t)x_2}{\alpha(t)+\beta(t) x_2+\gamma(t) x_3}\Big].
\end{aligned}
\end{cases}
\end{equation}

In this system, $x_i(t) (i=1,2,3)$ denote the  densities at time $t$ of three species, namely, PY1, PY2 and PR, respectively.
PY1 and  PY2 are two competitive  species and are  the food for the  predator species PR.
The set of equations \eqref{E1} describes the evolution of the system, where
\begin{itemize}
  \item At time $t$, $a_i(t) (i=1,2)$ are the intrinsic growth rates of  PY$i$, whereas  $a_3(t)$ is the death rate of PR.
  \item At time $t$, $b_{ij}(t) $  measures   the effect  PY$j$  has on the population of  PY$i$
  $ (i\ne j, i,j=1, 2)$,  meanwhile $b_{ii}(t)$ measures the inhibiting effect of the environment on PY$i$ $(i=1,2)$.
  \item An individual  predator consumes  amounts of prey, which are   Beddington-DeAngelis functional responses:
 $$ \frac{\sigma_1(t)x_1}{\alpha(t)+\beta(t) x_1+\gamma(t) x_3} \quad \text{ 
 and  } \quad  \frac{\sigma_2(t)x_2}{\alpha(t)+\beta(t) x_2+\gamma(t) x_3}$$
for some positive functions $\sigma_1, \sigma_2, \alpha, \beta$ and $\gamma$.  These amounts  contribute to the  growth  of the predator  a quantity of  the same form:
 $$\frac{\rho_1(t)x_1}{\alpha(t)+\beta(t) x_1+\gamma(t) x_3} + \quad  \frac{\rho_2(t)x_2}{\alpha(t)+\beta(t) x_2+\gamma(t) x_3}$$
for some positive functions $\rho_1$ and $\rho_2$.
\item All parameters are  continuous  and bounded above and below on $\mathbb R$ by some positive constants. 
\end{itemize}  
Problem  \eqref{E1} must be coupled with the initial condition for $x_i$:
$$x_i(0)=x_i^0 > 0, \hspace{1cm} i=1,2,3.$$

On one hand, we show survival conditions for all three species. When these species survive, a condition for the global stability of the system is presented. On the other hand, we show that when the density of both PY1 and PY2 is quite small, the predator species PR has not enough food and therefore falls into decay.

The organization of the  paper is  as follows.  Section \ref{preliminary} introduces some definitions such as permanence and stability for \eqref{E1}. 
  Section \ref{Sec3} shows the permanence of the three species and the decay of the predator species. We prove the existence of unique global positive solutions to  \eqref{E1} (Theorem \ref{Thm3.1}) and the existence of an invariant set under certain conditions (Theorem \ref{Thm3.2}).  Under these conditions, the system is permanent (Corollary  \ref{Cor3.1}). A condition for the decay of the predator species is also given (Theorem \ref{Thm3.3}). Section \ref{Sec4} investigates the global stability of solutions to \eqref{E1} (Theorem \ref{Thm4.1}). Finally, Section \ref{Sec5} gives some numerical examples.

\section{Definitions}  \label{preliminary}
Denote by  $\mathbb R^{\text 3}_+$ the positive cone of $\mathbb R^{\text 3}$, i.e.
$$\mathbb R^{\text 3}_+=\{ (x_1,x_2,x_3) \in \mathbb R^{\text 3}| x_i > 0, i=1,2,3\}.$$
Let $F$ be any function on $[0,\infty)$. For a brevity, sometime   we  write $F$ instead of  $F(t)$. In addition, denote
$$F^u= \sup_{0\leq t <\infty} F(t), \quad F^l=\inf_{0\leq t <\infty} F(t).$$ 

Let $x=(x_1, x_2, x_3)$  be a solution of \eqref{E1} with the  initial value $x^0=(x_1^0, x_2^0, x_3^0).$ 
\begin{definition} \label{defn1}
The system \eqref{E1} is said to be permanent if there exist $\delta_i >0  (i=1,2)$  such that 
$$ \delta_1\leq \liminf_{t \to \infty}x_i(t)\leq \limsup_{t \to \infty}x_i(t)\leq \delta_2 \hspace{1cm} $$ 
for all $x^0 \in \mathbb R^3_+$ and $i=1,2,3$.
\end{definition}
\begin{definition}  \label{defn2}
\begin{itemize}
  \item [{\rm (i)}] A set $A \subset \mathbb R^3_+$ is said to be positively invariant with respect to \eqref{E1} if for any $x^0\in A,$
     $$x(t)\in A, \hspace{2cm} 0\leq t<\infty.$$
  \item [{\rm (ii)}]
A set $A \subset \mathbb R^3_+$ is called  an ultimately bounded region of \eqref{E1} if for any solution $x$ of \eqref{E1}, there exists $T_1>0$ such that  
$$x(t) \in A, \hspace{2cm} T_1\leq t<\infty.$$
\end{itemize}
\end{definition}
\begin{definition}
A nonnegative solution $x^*$ of \eqref{E1} is called a global asymptotic stable solution  if it attracts any other solution $x$ of \eqref{E1} in the sense that
$$\lim_{t \to \infty}\sum_{i=1}^3|x_i(t)-x_i^*(t)|=0.$$
\end{definition}
\begin{remark}
It is easily seen that if  a solution of \eqref{E1}  is globally asymptotically  stable, then so are all solutions of \eqref{E1}. Thus, if there exists any  global asymptotic stable solution to \eqref{E1}, then we say that this system   is  globally asymptotically  stable.
\end{remark}

\section{Permanence and decay of species} \label{Sec3}
In this section, we  study  the permanence of the three species and the decay of the predator species. 
First, let us prove the existence of unique global positive solutions to  \eqref{E1}.

\begin{theorem}\label{Thm3.1}
For any  initial value $x^0\in \mathbb R_+^3$,  \eqref{E1} possesses  a unique global positive solution.
\end{theorem}

\begin{proof}
Since the functions in the right-hand sides of \eqref{E1} are local Lipschitz  continuous  on $\mathbb R_+^3$, there exists a unique local continuous solution $x=(x_1,x_2,x_3)$ defined on $[0,\tau)$, where either $\tau=\infty$  or $\tau$ is an explosion time, i.e. 
$$\tau<\infty \hspace{0.5cm} \text{and} \hspace{0.5cm} \lim_{t\to\tau} \|x(t)\|=\infty.$$

First, let us prove that $x(t)\in \mathbb R_+^3$ for  $0\leq t<\tau$. Denote by $\tau_0$ the first time in the period $[0,\tau)$ $x$ hits the boundary of $\mathbb R_+^3$, i.e.
\begin{equation}   \label{PT1}
\tau_0=\inf \{0\leq t<\tau| x_1(t)x_2(t)x_3(t)=0\}
\end{equation}
with a convention that $\inf \emptyset=\tau$. Because $(x_1(0),x_2(0),x_3(0))\in \mathbb R_+^3,$ we observe that
$$x_i(t)>0, \hspace{2cm} 0\leq t<\tau_0, \, i=1,2,3.$$
Hence, it is easily seen from  \eqref{E1} that
\begin{align}
x_1(t)=&x_1(0) \exp\Big\{\int_0^t \big[a_1(s)-b_{11}(s) x_1(s)-b_{12}(s)x_2(s) \notag\\
&- \frac{\sigma_1(s)x_3(s)}{\alpha(s)+\beta(s) x_1(s)+\gamma(s) x_3(s)}\big]ds\Big\}, \hspace{1cm} 0\leq t<\tau_0,  \label{PT2}\\
x_2(t)=&x_2(0) \exp\Big\{\int_0^t \big[a_2(s)-b_{21}(s) x_1(s)-b_{22}(s)x_2(s) \notag\\
&- \frac{\sigma_2(s)x_2(s)x_3(s)}{\alpha(s)+\beta(s) x_2(s)+\gamma(s) x_3(s)}\big]ds\Big\}, \hspace{1cm} 0\leq t<\tau_0,   \label{PT3}\\
x_3(t)=&x_3(0) \exp\Big\{\int_0^t \big[-a_3(s)+ \frac{\rho_1(s)x_1(s)}{\alpha(s)+\beta(s) x_1(s)+\gamma(s) x_3(s)}  \notag\\
&+\frac{\rho_2(s)x_2(s)}{\alpha(s)+\beta(s) x_2(s)+\gamma(s) x_3(s)}\big]ds\Big\}, \hspace{1cm} 0\leq t<\tau_0.  \label{PT4}
\end{align}
Since $x$ is continuous on $[0,\tau)$, \eqref{PT1}, \eqref{PT2}, \eqref{PT3} and \eqref{PT4}  imply that $\tau_0=\tau$. Thus, $x(t)\in \mathbb R_+^3$ for  $0\leq t<\tau$.

Second, let us prove that  $\tau=\infty$. Indeed, since $x(t)\in \mathbb R_+^3$ for  $0\leq t<\tau$,   \eqref{PT2}, \eqref{PT3} and \eqref{PT4} give
$$0<x_i(t) \leq x_i(0) \exp \Big\{\int_0^t a_i(s)ds\Big\}, \hspace{2cm} 0\leq t<\tau,i=1,2,$$
and
$$0<x_3(t) \leq x_3(0) \exp \Big\{\int_0^t \big[a_1(s)+\frac{\rho_1(s)}{\beta(s)}+\frac{\rho_2(s)}{\beta(s)}\big]ds\Big\}, \hspace{2cm} 0\leq t<\tau.$$
Because the functions $a_i (i=1,2,3), \rho_1, \rho_2$ and $\beta$ are continuous and bounded above and below on $\mathbb R$ by positive constants, we then obtain that 
$$\lim_{t\to \tau} \|x(t)\|<\infty.$$
By the definition of $\tau$, we conclude that $\tau=\infty$. This highlights that $x$ is a unique global positive solution of  \eqref{E1}
\end{proof}

Let us now show the existence of an invariant set and the permanence of all species.  
For $\epsilon\geq 0,$  put
\begin{equation*}
\begin{aligned}
& M_1^\epsilon=\frac{a_1^u}{b_{11}^l}+\epsilon, \ M_2^\epsilon=\frac{a_2^u}{b_{22}^l}+\epsilon,\\
&M_3^\epsilon=\frac{\rho_1^uM_1^\epsilon +\rho_2^uM_2^\epsilon -a_3^l\alpha^l}{a_3^l\gamma^l},\\
&m_1^\epsilon=\frac{(a_1^l-b_{12}^uM_2^0)(\alpha^l+\gamma^lM_3^0)-\sigma_1^uM_3^0}{b_{11}^u (\alpha^l+\gamma^lM_3^0)}-\epsilon,\\
&m_2^\epsilon=\frac{(a_2^l-b_{21}^uM_1^0)(\alpha^l+\gamma^lM_3^0)-\sigma_2^uM_3^0}{b_{22}^u (\alpha^l+\gamma^lM_3^0)}-\epsilon,\\
&m_3^\epsilon=\frac{(\rho_1^l-a_3^u\beta^u) m_1^0+(\rho_2^l-a_3^u\beta^u)m_2^0-2a_3^u\alpha^u}{2a_3^u\gamma^u}-\epsilon,
\end{aligned}
\end{equation*}
and
\begin{equation}\label{E2}
\Gamma=\{(x_1, x_2, x_3) \in \mathbb R^3| \ m_i^0 \leq  x_i  \leq  M_i^0, \ i=1,2,3\}.
\end{equation}

\begin{theorem} \label{Thm3.2}
Assume that $M_3^0> 0$ and $ m_i^0>0\, (i=1,2,3).$ 
Then,  the set $\Gamma$ defined by \eqref{E2}
is positively invariant with respect to  \eqref{E1}.
\end{theorem}
\begin{proof}
Since $M_3^0> 0, m_i^0>0,$ it is easily seen that  $M_i^0>m_i^0>0$ $ (i=1,2,3)$. Thereby, $\Gamma$ is a non-empty set in  the positive cone of $\mathbb R^3$.

Let $x^0\in\Gamma$. It suffices to prove that  $x(t)\in \Gamma$ for   $0\leq t <\infty$. For this, we use  a  fact that the  equation
\begin{equation*}
\begin{cases}
X'(t)=A(t) X(t) [B-X(t)] \hspace{1cm} (B \neq 0),\\
X(0)=X^0
\end{cases}
\end{equation*}
has the explicit solution
$$X(t) = \frac{BX^0  e^{\int_{0}^t A(s)Bds}}{X^0[e^{\int_{0}^t A(s)Bds}-1 ]+B}.$$

First, let us  give upper estimates for $x_i (i=1,2,3)$. The first equation of  \eqref{E1} gives 
\begin{equation*}
\begin{aligned}
x_1'(t) &\leq x_1(t)[a_1(t)-b_{11}(t)x_1(t)]\\
& \leq x_1(t)[a_1^u-b_{11}^lx_1(t)]\\
&=b_{11}^lx_1(t)(M_1^0-x_1).
\end{aligned}
\end{equation*}
The comparison theorem and the inequality $0<x_1^0\leq M_1^0$ then provide that  
\begin{align}
x_1(t) &\leq \frac{x_1^0 M_1^0 e^{a_{1}^ut}}{x_1^0(e^{a_1^ut}-1) +M_1^0} \notag\\
&\leq \frac{x_1^0M_1^0 e^{a_1^ut}}{x_1^0(e^{a_1^ut}-1)+M_1^0} \notag\\
&\leq  M_1^0, \hspace{2cm}0\leq  t <\infty. \label{E3}
\end{align}

Similarly,  
\begin{equation}\label{E4}
x_2(t)\leq M_2^0, \hspace{2cm}  0\leq  t <\infty.
\end{equation}

Substituting   \eqref{E3} and  \eqref{E4} into the third equation of  \eqref{E1}, we obtain  that
\begin{align}
x_3'&	\leq -a_3^lx_3 + \frac{\rho_1^ux_1x_3}{\alpha^l+\beta^l x_1+\gamma^l x_3}+\frac{\rho_2^ux_2x_3}{\alpha^l+\beta^l x_2+\gamma^l x_3}\notag \\
&	\leq -a_3^lx_3 + \frac{(\rho_1^uM_1^0 +\rho_2^uM_2^0 )x_3}{\alpha^l+\gamma^l x_3}\notag\\
&	= \frac{x_3 \left[(\rho_1^uM_1^0 +\rho_2^uM_2^0 -a_3^l\alpha^l)-a_3^l\gamma^l x_3\right]}{\alpha^l+\gamma^l x_3}\notag\\
&=\frac{a_3^l\gamma^l}{\alpha^l+\gamma^l x_3} x_3(M_3^0- x_3). \label{E5}
\end{align}
The comparison theorem again provides that  
\begin{equation}\label{E6}
x_3(t) \leq  \frac{x_3^0 M_3^0 e^{M_3^0\int_{0}^tC_1(s)ds}}{x_3^0\Big[e^{M_3^0\int_{0}^tC_1(s)ds}-1\Big]+M_3^0} \leq M_3^0, \hspace{1cm}  0\leq  t <\infty
\end{equation}
(note that  $0<x_3^0 \leq M_3^\epsilon$), where 
\begin{equation} \label{E56}
C_1(t)=\frac{a_3^l\gamma^l}{\alpha^l+\gamma^l x_3(t)}.
\end{equation}

Second, let us   give  lower estimates for $x_i (i=1,2,3)$.  Thanks  to  \eqref{E1}, \eqref{E3}, \eqref {E4} and \eqref{E6},  
\begin{equation*}
\begin{aligned}
x_1'(t) \geq &x_1 \Big(a_1^l-b_{11}^ux_1-b_{12}^uM_2^0 -\frac{\sigma_1^uM_3^0 }{\alpha^l+\beta^l x_1+\gamma^lM_3^0}\Big)\\
\geq &x_1 \Big(a_1^l-b_{11}^ux_1-b_{12}^uM_2^0 -\frac{\sigma_1^uM_3^0 }{\alpha^l+\gamma^lM_3^0}\Big)\\
 \geq &x_1 \Big[\frac{(a_1^l-b_{12}^uM_2^0)(\alpha^l+\gamma^lM_3^0)-\sigma_1^uM_3^0}{\alpha^l+\gamma^lM_3^0}-b_{11}^u x_1\Big]\\
=&b_{11}^u x_1(m_1^0 -x_1).\\
\end{aligned}
\end{equation*}
 The comparison theorem then gives  
\begin{equation}\label{E7}
x_1(t) \geq  m_1^0,  \hspace{2cm}  0\leq t <\infty.
\end{equation}

Similarly, 
$$ x_2(t) \geq m_2^0, \hspace{2cm}  0\leq t <\infty.
$$

In the meantime,  the last equation of \eqref{E1}  gives
\begin{equation*}
\begin{aligned}
x_3'
=&-a_3x_3+ \frac{\rho_1x_1x_3}{\alpha+\beta x_1+\gamma x_3}+\frac{\rho_2x_2x_3}{\alpha+\beta x_2+\gamma x_3}\\
\geq &-a_3^u x_3+\frac{\rho_1^l m_1^0 x_3}{\alpha^u+\beta^u m_1^0+\gamma^u x_3}+\frac{\rho_2^l m_2^0 x_3}{\alpha^u+\beta^u m_2^0+\gamma^u x_3}\\
\geq & x_3\left[-a_3^u+ \frac{\rho_1^l m_1^0+\rho_2^lm_2^0}{2\alpha^u+\beta^u (m_1^0+m_2^0)+2\gamma^u x_3}\right]\\
= &\frac{x_3}{2\alpha^u+\beta^u (m_1^0+m_2^0)+2\gamma^u x_3}\\
&\times [\rho_1^l m_1^0+\rho_2^lm_2^0-a_3^u\{2\alpha^u+\beta^u(m_1^0+m_2^0)\} -2a_3^u\gamma^u x_3 ] \\
=&\frac{2a_3^u\gamma^u}{2\alpha^u+\beta^u (m_1^0+m_2^0)+2\gamma^u x_3} x_3 (m_3^0-x_3),
\end{aligned}
\end{equation*}
here we used the inequality: 
$$\frac{x_1}{y_1}+\frac{x_2}{y_2} \geq \frac{x_1+x_2}{y_1+y_2}, \hspace{1cm} x_i, y_i>0, i=1,2.$$
The comparison theorem then provides that 
$$x_3(t) \geq m_3^0, \hspace{2cm}   0\leq t <\infty.$$ 

By the above arguments, we have thus concluded that $x(t)\in \Gamma$ for  $0\leq t<\infty$. The proof is now complete.
\end{proof}

\begin{corollary}[Permanence of species]\label{Cor3.1} 
Assume that $M_3^0> 0$ and $ m_i^0>0\, (i=1,2,3).$  Then, the solution $x=(x_1,x_2,x_3)$ of \eqref{E1} with $x(0)=x^0\in \mathbb R_+^3$ satisfies
$$ m_i^0 \leq \liminf_{t \to \infty}x_i(t)\leq  \limsup_{t \to \infty}x_i(t)\leq M_i^0, \hspace{1cm} i=1,2,3.$$
This means that the system \eqref{E1} is permanent and  $\Gamma$  is an ultimately bounded region.
\end{corollary}

\begin{proof}
According to the proof for  Theorem \ref {Thm3.2}, we have 
$$x_i(t) \leq \frac{x_i^0 M_i^0 e^{a_i^u t }}{x_i^0[e^{a_i^u t }-1] +M_i^0}, \hspace{1cm}  i=1, 2.$$
Thereby, 
$$\limsup_{t \to \infty}x_i(t)\leq M_i^0, \hspace{1cm}  i=1, 2.$$
As a consequence,  for any $\epsilon >0,$ there exists $0\leq t_1<\infty$ such that for $i=1,2,$
$$x_i(t)< M_i^\epsilon, \hspace{2cm}    t_1\leq t<\infty.$$

In the meantime, using the  same argument as in   \eqref{E5} and \eqref{E6}, we observe that 
\begin{equation}\label{E8}
x_3(t) \leq  \frac{M_3^\epsilon x_3^1 e^{M_3^\epsilon\int_{t_1}^tC_1(s)ds}}{x_3^1\big[e^{M_3^\epsilon\int_{t_1}^tC_1(s)ds}-1\big]+M_3^\epsilon}, \hspace{2cm}    t_1\leq t<\infty,
\end{equation}
  where $x_3^1=x_3(t_1)$ and $C_1$ is defined by \eqref{E56}.  This implies that
$$x_3(t) \leq \max\{M_3^\epsilon, x_3^1\}, \hspace{2cm}    t_1\leq t<\infty.$$
Clearly, by  \eqref{E56}, 
  $$\inf_{t\geq t_1} C_1(s) >0.$$ 
 Taking the limit  as $t\to \infty$ in   \eqref{E8}, it therefore follows that 
 $$\limsup_{t \to \infty} x_3(t) \leq M_3^\epsilon.$$
Again taking the limit  as $\epsilon \to 0,$  we arrive at 
$$\limsup_{t \to \infty} x_3(t) \leq M_3^0.$$

   Similarly, we obtain that 
$$\liminf_{t \to \infty}x_i(t)\geq m_i^0, \hspace{1cm}  i=1,2,3.$$
 This completes the proof.
\end{proof}

Finally, let us  give a sufficient condition for the decay of the predator species.
\begin{theorem}[Decay of predator species] \label{Thm3.3}
Assume that $M_3^0<0.$ Let  $x=(x_1,x_2,x_3)$ be the solution of  \eqref{E1} with $x(0)=x^0\in \mathbb R_+^3$. Then, 
$$\lim\limits_{t \to \infty}x_3(t)=0.$$
 This means that  the predator species falls into decay.
\end{theorem}

\begin{proof}
In the proof for  Corollary  \ref{Cor3.1},  we already proved that for any $\epsilon>0$, there exists $0\leq t_1<\infty$ such that for $i=1,2$, 
$$x_i(t)< M_i^\epsilon, \hspace{2cm}    t_1\leq t<\infty.$$

Meanwhile, 
using the same arguments as in  \eqref{E5}, we obtain that 
\begin{equation} \label{E9}
x_3'(t)\leq \frac{a_3^l\gamma^l}{\alpha^l+\gamma^l x_3} x_3(t) [M_3^\epsilon- x_3(t)], \hspace{1cm} t_1 \leq t<\infty.
\end{equation}
Since  $M_3^0<0$, we choose $\epsilon>0$ such that $M_3^\epsilon<0$. Therefore, \eqref{E9} gives
$$x_3'(t)<0, \hspace{1cm} t_1 \leq t<\infty.$$
The function  $x_3$ is hence  decreasing  in time.  As a consequence, there exists $C\geq 0$ such that 
$$\lim\limits_{t \to \infty}x_3(t)=C \quad \text{ and } \quad  C \leq x_3(t) \leq x_3(t_1), \hspace{1cm}  t_1 \leq t<\infty.$$
If $C>0,$ then by \eqref{E9} there would exist $\mu>0$ such that $x_3'(t)<-\mu$ for $ t_1 \leq t<\infty$.   
As a consequence,  
$$x_3(t)<-\mu (t-t_1) +x_3(t_1), \hspace{2cm}  t_1 \leq t<\infty.$$
 This leads to a contradiction: 
 $$0\leq \lim\limits_{t \to \infty}x_3(t)=-\infty.$$
Therefore, $C=0$ and thus $\lim\limits_{t \to \infty}x_3(t)=0.$
\end{proof}

\begin{remark}  \label{rm1}
Suppose that all the  parameters in the system \eqref{E1} are independent of time, i.e. they are constants.

When the second prey species PY2 and the predator species PR are absent from the system, the density $x_1$ of the first prey species PY1  is governed by the logistic equation
$$x_1'=x_1(a_1-b_{11}x_1).$$
It is easily seen that for this case  $x_1$   tends towards the  equilibrium state $x_1^*=\frac{a_1}{b_{11}}$. Similarly, when the first prey species and the predator are absent from the system, the density $x_2$ of the second prey species tends towards the  equilibrium state $x_2^*=\frac{a_2}{b_{22}}$. 

As shown in Corollary \ref{Cor3.1}, the solution $x=(x_1,x_2,x_3)$ of  \eqref{E1} with $x(0)=x^0\in \mathbb R_+^3$
satisfies
$$\limsup_{t\to \infty}x_i(t) \leq x_i^*,  \hspace{1cm} i=1,2.$$
The condition $M_3^0<0$ in Theorem \ref{Thm3.3} is therefore   equivalent  to the condition: 
$$\frac{\rho_1 x_1^*+\rho_2x_2^*}{\alpha}<a_3.$$
This condition has  biological meaning. When the dead rate $a_3$ of the predator species  exceeds  the limit of the contribution of two prey species to the growth of the predator species, the predator species falls into decay. The large dead rate may result from epidemic disease or random factors (see \cite{Linh,TVT1,TVT3}).
 \end{remark}

\section{Global stability} \label{Sec4}
In this section, we  show the  global asymptotic stability of solutions to  \eqref{E1}. 
The following lemma is useful.
\begin{lemma}\label{lem2}
Let $m$ be a real number and $F$ be a nonnegative function defined on $[m, \infty)$ such that $F$ is integrable and uniformly continuous on $[m, \infty).$ Then,
$$\lim_{t \to \infty}F(t)=0.$$
\end{lemma}
\begin{proof}
Suppose the contrary that $F(t) \not\rightarrow 0$ as $t \rightarrow \infty$. Given $\epsilon>0$. There would then exist a sequence $\{t_{n}\}_{n=1}^\infty, t_n\geq a$ such that $t_{n} \rightarrow \infty$ as $n \rightarrow \infty$ and $F(t_{n}) \geq \varepsilon$ for  $n=1,2,3,\dots$ 
By the uniform continuity of $F$, there exists  $\delta > 0$ such that 
$$ |F(t_{n}) - F(t)| \leq \frac{\varepsilon}{2},  \hspace{2cm} t_{n} \leq t \leq t_{n}+\delta, n =1,2,3,\dots$$
 Thus,  
\begin{equation*}
\begin{aligned}
 F(t)=&|F(t_{n}) - [F(t_{n})- F(t)]| \\
 \geq &|F(t_{n})| - |F(t_{n})- F(t)| 
\geq \varepsilon - \frac{\varepsilon}{2} = \frac{\varepsilon}{2}.
 \end{aligned}
\end{equation*}
Therefore,  for each $n =1,2,3,\dots,$
$$ \int_{t_{n}}^{t_{n}+\delta} F(t) dt = \int_{t_{n}}^{t_{n}+\delta} F(t) dt \geq \frac{\varepsilon \delta}{2} > 0. $$

 On the other hand, since the Riemann integral $\int_{a}^{\infty} F(t) dt$ exists, the integral $ \int_{t_{n}}^{t_{n}+\delta} F(t) dt$  converges to 0 as $n \rightarrow \infty.$ We then arrive at a  contradiction:
$$0\geq  \frac{\varepsilon \delta}{2}.$$
Thus, $\lim_{t\to\infty} F(t)=0$.
\end{proof}

\begin{theorem} [Global stability]\label{Thm4.1}
Assume that $M_3^0> 0$ and $ m_i^0>0\, (i=1,2,3).$ Let $x^*$ be a solution of \eqref{E1} with $x^*(0) \in \mathbb R_+^3$. Let $\epsilon>0$ be  sufficiently small  such that  $m_i^\epsilon >0\, (i=1,2,3)$  and 
\begin{equation} \label{E11}
\begin{cases}
\begin{aligned}
\sup_{t \geq 0}&\Big\{b_{21}+\frac{\alpha \rho_1+(\gamma \rho_1+\beta \sigma_1)M_3^\epsilon}{u_1(m_1^\epsilon, M_3^\epsilon)}-b_{11}\Big\} <0,\\
\sup_{t \geq 0}&\Big\{b_{12}+\frac{\alpha \rho_2+(\gamma \rho_2+\beta \sigma_2)M_3^\epsilon}{u_2(m_2^\epsilon, M_3^\epsilon)}-b_{22}\Big\}<0,\\
\sup_{t \geq 0}&\Big\{\frac{\sigma_1(\alpha+\beta M_1^\epsilon)}{u_1(M_1^\epsilon,m_3^\epsilon)}+\frac{\sigma_2(\alpha+\beta M_2^\epsilon)}{u_2(M_2^\epsilon,m_3^\epsilon)}-\frac{\gamma \rho_1 m_1^\epsilon}{u_1(m_1^\epsilon,M_3^\epsilon)}-\frac{\gamma \rho_2 m_2^\epsilon}{u_2(m_2^\epsilon,M_3^\epsilon)}\Big\}<0,
\end{aligned}
\end{cases}
\end{equation}
where $u_i(u,v)=(\alpha+\beta x_i^{*}+\gamma x_3^*)(\alpha+\beta u+\gamma v) \, (i=1,2).$
Then, $x^*$ is globally asymptotically  stable.
\end{theorem}

\begin{proof}
Let $x$ be any other solution of \eqref{E1} with $x(0) \in \mathbb R_+^3$. Thanks to  Corollary \ref{Cor3.1}, there exists $T_1>0$ such that for $i=1,2,3$,
$$m_i^\epsilon<x_i(t)<M_i^\epsilon, \hspace{2cm} T_1\leq t<\infty.$$

Consider a Lyapunov function defined by
 $$V(t)=\sum_{i=1}^3 |\ln x_i-\ln x_i^*|, \hspace{2cm} 0\leq t<\infty.$$
A direct calculation of the right derivative $D^{+}V$ of  $V$ along the solutions of \eqref{E1} gives
\begin{align*}
D^{+}V(t)	
=&\sum_{i=1}^3 \mathop{\rm sgn}(x_i-x_i^*) \left( \frac{x_i'}{x_i}-\frac{{x_i^*}'}{x_i^*}\right)
\\
=&\mathop{\rm sgn}(x_1-x_1^*)\Big[ -b_{11}(x_1-x_1^*)-b_{12}(x_2-x_2^*)\\
&-\sigma_1\Big(\frac{x_3}{\alpha+\beta x_1+\gamma x_3}-\frac{x_3^*}{\alpha+\beta x_1^*+\gamma x_3^*}\Big)\Big] \\
&+\mathop{\rm sgn}(x_2-x_2^*)\Big[ -b_{21}(x_1-x_1^*)-b_{22}(x_2-x_2^*)\\
&-\sigma_2\Big(\frac{x_3}{\alpha+\beta x_2+\gamma x_3}-\frac{x_3^*}{\alpha+\beta x_2^*+\gamma x_3^*}\Big)\Big] \\
&+\mathop{\rm sgn}(x_3-x_3^*) \Big(\frac{\rho_1 x_1}{\alpha+\beta x_1+\gamma x_3}-\frac{\rho_1 x_1^*}{\alpha+\beta x_1^*+\gamma x_3^*}\\
&+\frac{\rho_2 x_2}{\alpha+\beta x_2+\gamma x_3}-\frac{\rho_2 x_2^*}{\alpha+\beta x_2^*+\gamma x_3^*} \Big)\\
\leq&(b_{21}-b_{11})|x_1-x_1^*|+(b_{12}-b_{22})|x_2-x_2^*|\\
&- \sigma_1 \mathop{\rm sgn} (x_1-x_1^*) \frac{\alpha (x_3-x_3^*)+\beta(x_1^*x_3-x_1x_3^*)}{u_1(x_1,x_3)}\\
&- \sigma_2 \mathop{\rm sgn} (x_2-x_2^*) \frac{\alpha (x_3-x_3^*)+\beta(x_2^*x_3-x_2x_3^*)}{u_2(x_2,x_3)}\\
&+\mathop{\rm sgn}(x_3-x_3^*) \frac{\alpha \rho_1(x_1-x_1^*)+\gamma \rho_1 (x_1 x_3^*-x_1^* x_3)}{u_1(x_1, x_3)}\\
&+\mathop{\rm sgn}(x_3-x_3^*) \frac{\alpha \rho_2(x_2-x_2^*)+\gamma \rho_2 (x_2 x_3^*-x_2^* x_3)}{u_2(x_2, x_3)}.
\end{align*}
Using the expression 
$$x_ix_3^*-x_i^*x_3=x_i(x_3^*-x_3)+x_3(x_i-x_i^*)\quad (i=1,2),$$ 
we observe that
\begin{align}\label{E12}
D^+V(t)
 \leq&(b_{21}-b_{11})|x_1-x_1^*|+(b_{12}-b_{22})|x_2-x_2^*| 										\notag \\
&- \sigma_1 \mathop{\rm sgn} (x_1-x_1^*) \frac{(\alpha+\beta x_1) (x_3-x_3^*)-\beta x_3(x_1-x_1^*)}{u_1(x_1,x_3)}					\notag\\
&- \sigma_2 \mathop{\rm sgn} (x_2-x_2^*) \frac{(\alpha+\beta x_2) (x_3-x_3^*)-\beta x_3(x_2-x_2^*)}{u_2(x_2,x_3)}					\notag\\
&+\mathop{\rm sgn}(x_3-x_3^*) \frac{(\alpha +\gamma x_3)\rho_1(x_1-x_1^*)-\gamma \rho_1 x_1(x_3- x_3^*)}{u_1(x_1, x_3)}					\notag\\
&+\mathop{\rm sgn}(x_3-x_3^*)  \frac{(\alpha +\gamma x_3)\rho_2(x_2-x_2^*)-\gamma \rho_2 x_2(x_3- x_3^*)}{u_2(x_2, x_3)}						\notag\\
\leq&\left[b_{21}+\frac{(\alpha +\gamma M_3^\epsilon)\rho_1+\beta \sigma_1M_3^\epsilon}{u_1(m_1^\epsilon, M_3^\epsilon)}-b_{11}\right]|x_1-x_1^*|		\notag\\
&+\left[b_{12}+\frac{(\alpha +\gamma M_3^\epsilon)\rho_2+\beta \sigma_2M_3^\epsilon}{u_2(m_2^\epsilon, M_3^\epsilon)}-b_{22}\right]|x_2-x_2^*|		\notag\\
&+\Big[\frac{\sigma_1(\alpha+\beta M_1^\epsilon)}{u_1(M_1^\epsilon,m_3^\epsilon)}+\frac{\sigma_2(\alpha+\beta M_2^\epsilon)}{u_2(M_2^\epsilon,m_3^\epsilon)}-\frac{\gamma \rho_1 m_1^\epsilon}{u_1(m_1^\epsilon,M_3^\epsilon)} -\frac{\gamma \rho_2 m_2^\epsilon}{u_2(m_2^\epsilon,M_3^\epsilon)} \Big]   \notag\\
&\times|x_3-x_3^*|                                                               \notag\\
=&\left[b_{21}+\frac{\alpha \rho_1+(\gamma \rho_1+\beta \sigma_1)M_3^\epsilon}{u_1(m_1^\epsilon, M_3^\epsilon)}-b_{11}\right]|x_1-x_1^*|                    \notag\\
&+\left[b_{12}+\frac{\alpha \rho_2+(\gamma \rho_2+\beta \sigma_2)M_3^\epsilon}{u_2(m_2^\epsilon, M_3^\epsilon)}-b_{22}\right]|x_2-x_2^*|                     \notag\\
&+\Big[\frac{\sigma_1(\alpha+\beta M_1^\epsilon)}{u_1(M_1^\epsilon,m_3^\epsilon)}+\frac{\sigma_2(\alpha+\beta M_2^\epsilon)}{u_2(M_2^\epsilon,m_3^\epsilon)}-\frac{\gamma \rho_1 m_1^\epsilon}{u_1(m_1^\epsilon,M_3^\epsilon)}                       -\frac{\gamma \rho_2 m_2^\epsilon}{u_2(m_2^\epsilon,M_3^\epsilon)} \Big]\notag\\
&\times |x_3-x_3^*|, \hspace{2cm}  T_1\leq t<\infty.
\end{align}

Combining \eqref{E11} and \eqref{E12}, it is easily seen that there exists $\mu>0$ such that 
\begin{equation}\label{E13}
D^{+}V(t) \leq -\mu \sum_{i=1}^3 |x_i-x_i^*|, \hspace{2cm}  T_1\leq t<\infty.
\end{equation}
Integrating both the hand sides of \eqref{E13} from  $T_1$ to $t\geq T_1$, it follows that 
$$V(t)+\mu \int_{T_1}^{t}\sum_{i=1}^3 |x_i-x_i^*|ds \leq V(T_1) <\infty.$$
Hence, 
$$ \int_{T_1}^{t}\sum_{i=1}^3 |x_i-x_i^*|ds \leq \mu^{-1}V(T_1) < \infty, \hspace{1cm} T_1 \leq t<\infty.$$
Thus,  
$ \sum_{i=1}^3 |x_i-x_i^*| \in L^1([T_1, \infty)).$

On the other hand, by the ultimate boundedness of  $x_i^*$ and $x_i, $  the functions  $x_i$ and $ x_i^* \, (i= 1,2,3)$ have bounded derivatives at time $t \geq T_1$. As a consequence, 
$ \sum_{i=1}^3 |x_i-x_i^*|$ is uniformly continuous on $[T_1, \infty)$. Due to Lemma \ref{lem2}, we conclude that  
$$\underset{t \to \infty}{\lim}\sum_{i=1}^3 |x_i(t)-x_i^*(t)|=0.$$
The solution $x^*$ is thus globally asymptotically  stable.
\end{proof}

\begin{remark}  \label{rm2}
Consider the case where all parameters in  \eqref{E1} are constants. Since 
$$u_i(u,v)=(\alpha+\beta x_i^{*}+\gamma x_3^*)(\alpha+\beta u+\gamma v)\geq \alpha (\alpha+\beta u+\gamma v),$$
where $u_i (i=1,2) $ are defined in  Theorem \ref{Thm4.1}, it is easily seen that the condition \eqref{E11} can be replaced by a simple one:
\begin{equation} \label{E30}
\begin{cases}
\begin{aligned}
&b_{21}+\frac{\alpha \rho_1+(\gamma \rho_1+\beta \sigma_1)M_3^0}{\alpha (\alpha+\beta m_1^0+\gamma M_3^0)}-b_{11} <0,\\
&b_{12}+\frac{\alpha \rho_2+(\gamma \rho_2+\beta \sigma_2)M_3^0}{\alpha (\alpha+\beta m_2^0+\gamma M_3^0) }-b_{22}<0,\\
&\frac{\sigma_1(\alpha+\beta M_1^0)}{\alpha (\alpha+\beta M_1^0+\gamma m_3^0)}+\frac{\sigma_2(\alpha+\beta M_2^0)}{\alpha (\alpha+\beta M_2^0+\gamma m_3^0) }\\
&-\frac{\gamma \rho_1 m_1^0}{\alpha (\alpha+\beta m_1^0+\gamma M_3^0)  }-\frac{\gamma \rho_2 m_2^0}{\alpha (\alpha+\beta m_2^0+\gamma M_3^0)  }<0.
\end{aligned}
\end{cases}
\end{equation}
In applications, the condition \eqref{E30} can be easily checked.

 \end{remark}

\section{Numerical examples} \label{Sec5}
In this section, we present  some numerical examples. First, we give examples
which show the permanence of all species and the  decay of the predator species; second, examples which suggest possibility of global stability.

\subsection{Examples for permanence and decay}
Let us first observe an example that show that, if the dead rate of the predator species is sufficiently
small, then all the species grow stably.

Set $a_1= a_2=8, a_3=0.5, b_{11}= b_{22}=4,  b_{12}=4-\frac{1}{1+t^2}, b_{21}=3+\frac{1}{1+t^2}, \sigma_1=4-\frac{1}{1+t}$ $ \sigma_2=3+\frac{1}{1+t}, \alpha=1+\frac{1}{1+t}, 
\beta=\frac{1+t}{2+t}, 
\gamma=1,
\rho_1= \rho_2=2$. 

It is easily seen that the assumptions of Theorem \ref{Thm3.2} and Corollary \ref{Cor3.1}  are satisfied. 
We calculate 5000 trajectory values of $(x_1,x_2,x_3)$ at $T = 100$, where  all  initial values are randomly generated
in the rectangular parallelepiped  $[0.2,10.2]\times[0.2,10.2] \times [0.2,10.2]$.  

\begin{figure}[H] 
\begin{center}
\includegraphics[scale=0.7]{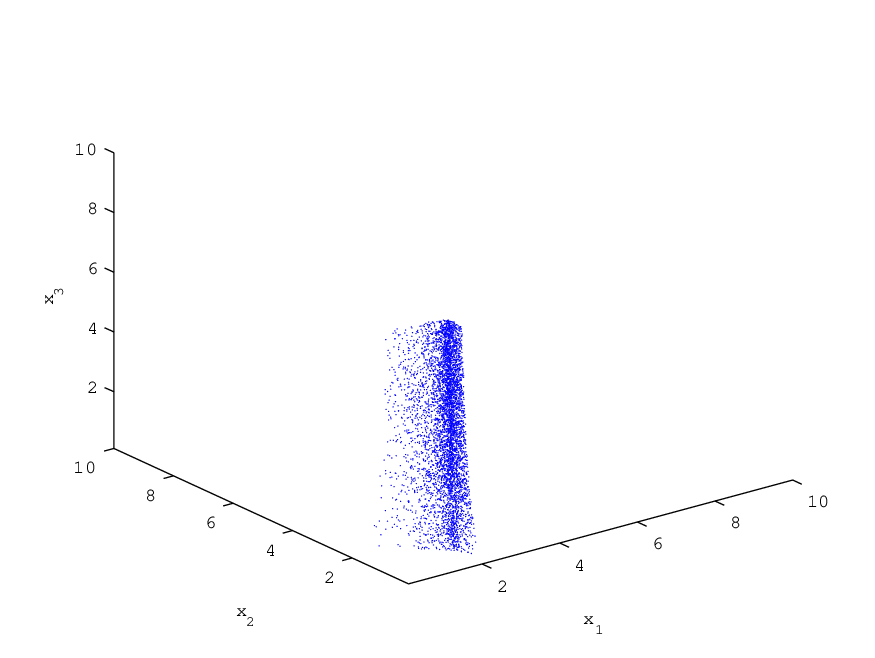}
\end{center}
\vspace*{9pt}
\caption{Ultimately bounded region of \eqref{E1}. Corresponding to 5000 different initial values, the densities of three species ultimately concentrate on this region.}
\label{Fig1}
\end{figure}

Figure \ref{Fig1} shows an ultimately bounded region of \eqref{E1}. It suggests that the system is permanent. 

Let us next observe an example  showing that, if the dead rate of the predator species is sufficiently
large, then the predator species falls into the decay.

Set $a_1=1, a_2=2, a_3=2,   
\rho_1=2, \rho_2=1$.       
   The other parameters of \eqref{E1} are taken the same to the above example. 

\begin{figure}[H] 
\begin{center}
\includegraphics[height=5.8cm,width=10cm]{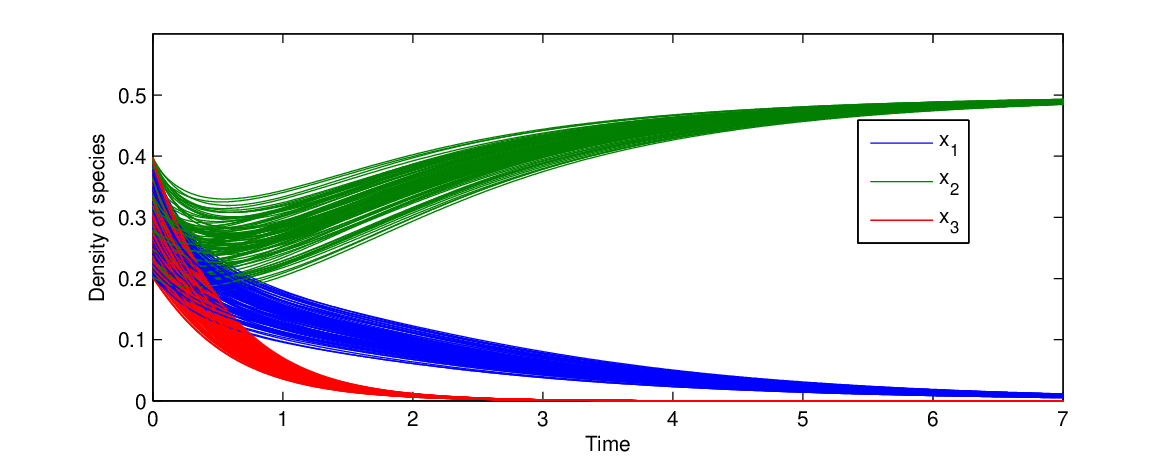}
\end{center}
\caption{The density  of three species. All 100 trajectories of $x_3$ approach zero in time, i.e. the predator species falls into decay.}
\label{Fig2}
\end{figure}

Similarly to the above example, these parameters satisfy the assumptions of Theorem \ref{Thm3.3}. We calculate 100 trajectories of $(x_1,x_2,x_3)$ on  $[0,30]$, where  all  initial values are randomly generated
in the rectangular parallelepiped  $[0.2,0.4]\times[0.2,0.4] \times [0.2,0.4]$.  
 Figure \ref{Fig2} shows that all the trajectories of the density $x_3$ of the predator species tend to $0$ in time.

\subsection{Example for global stability}
Let us  observe an example that show the global stability of solutions of \eqref{E1}. 
Set $a_1=4, a_2=10, a_3=0.5, b_{11}= b_{22}=4, \alpha=\beta=\gamma=  b_{12}= b_{21}=1; \sigma_1=\sigma_2=0.1,
\rho_1= \rho_2=2$. A set of 150 initial values of $(x_1(0),x_2(0),x_3(0))$ is randomly generated in the rectangular parallelepiped $[0.1,20.1]\times[0.1,20.1] \times [0.1,20.1]$.

The system \eqref{E1} is hence of autonomous differential equations. It is easily checked that  \eqref{E1} possesses  a unique stationary solution $x^*=(x_1^*,x_2^*,x_3^*)$
and that the assumptions of Theorem \ref{Thm4.1} is satisfied (note that the condition \eqref{E11} is replaced by \eqref{E30}).

\begin{figure}[H] 
\begin{center}
\includegraphics[height=5.8cm,width=10cm]{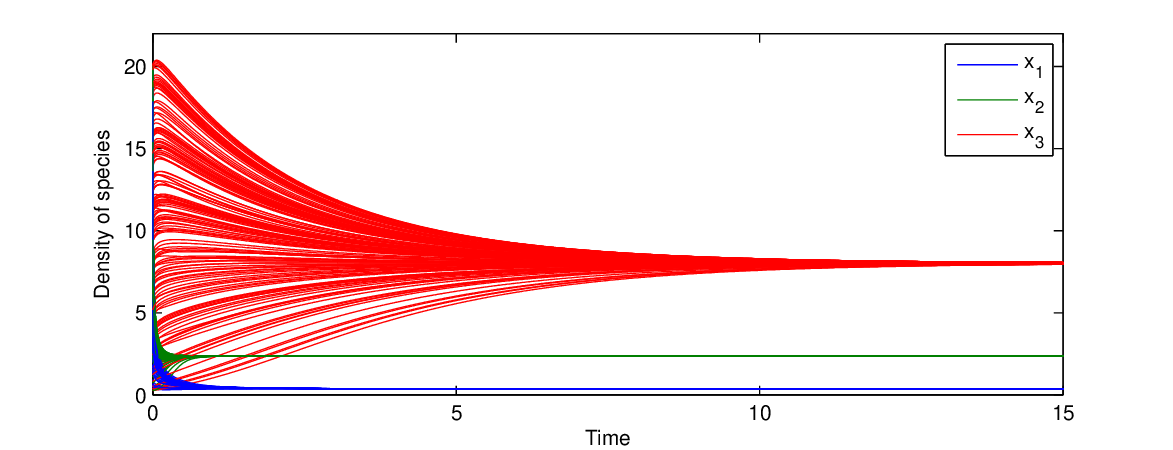}
\end{center}
\caption{Global stability of \eqref{E1}. All the 150 solutions of \eqref{E1} tend to a fixed point, i.e. the stationary solution $x^*$.}
\label{Fig3}
\end{figure}

 We plot 150 trajectories of   $(x_1,x_2,x_3)$ corresponding to the set of initial values  on $[0,15]$ in     
Figure \ref{Fig3}. It is seen that  all the solutions of \eqref{E1} converge to $x^*$ in time. 

\section{Conclusions}
The non-periodic version of the one-predator two-prey system model presented in 
\cite{SICE} has been studied. We showed the variation of density of each species caused other species by giving  conditions for permanence, stability and decay of species. 

Our system may be applied to the study of many biological  phenomena  in the real world. For examples, 
it can be used to predict the growth of lion, buffalo, and  zebra species in the Serengeti ecosystem in Africa, or to propose a plan for the forest exploitation on deciduous and coniferous trees (e.g., the mixed coniferous deciduous forest of white birch aspens fir and spruce trees in the Superior National Forest,  Minnesota, United States (\cite{Minnesota})).



\end{document}